\documentclass[12pt]{article} 
\usepackage{amsfonts,amsmath,latexsym,amssymb,mathrsfs,amsthm,comment}
\usepackage{slashbox}
\usepackage{caption}

\let\OLDthebibliography\thebibliography
\renewcommand\thebibliography[1]{
  \OLDthebibliography{#1}
  \setlength{\parskip}{3pt}
  \setlength{\itemsep}{0pt plus 0.3ex}
}

\def\numberlikeadb{\global\def\theequation{\thesection.\arabic{equation}}}
\numberlikeadb
\newtheorem{theorem}{Theorem}[section]

\newtheorem{corollary}[theorem]{Corollary}

\newtheorem{remark}[theorem]{Remark}

\usepackage{color}

\usepackage{lscape}
\usepackage{caption}
\usepackage{multirow}
\allowdisplaybreaks
\begin{document}

\title{On the moments of the variance-gamma distribution
}
\author{Robert E. Gaunt\footnote{Department of Mathematics, The University of Manchester, Oxford Road, Manchester M13 9PL, UK}  
}

\date{} 
\maketitle

\vspace{-10mm}

\begin{abstract}We obtain new closed-form formulas for the moments and absolute moments of the variance-gamma distribution. We thus deduce new formulas for the moments and absolute moments of the product of two correlated zero mean normal random variables.

\end{abstract}


\noindent{{{\bf{AMS 2010 Subject Classification:}}} Primary 60E05; 62E15

\section{Introduction}

The variance-gamma (VG) distribution with parameters $\nu > -1/2$, $0\leq|\beta|<\alpha$, $\mu \in \mathbb{R}$, which we denote by $\mathrm{VG}(\nu,\alpha,\beta,\mu)$, has probability density function (PDF)
\begin{equation}\label{vgpdf} p(x) = M \mathrm{e}^{\beta (x-\mu)}|x-\mu|^{\nu}K_{\nu}(\alpha|x-\mu|), \quad x\in \mathbb{R},
\end{equation}
where 
\[M=M_{\nu,\alpha,\beta}=\frac{(\alpha^2-\beta^2)^{\nu+1/2}}{\sqrt{\pi}(2\alpha)^\nu \Gamma(\nu+1/2)}\]
 is the normalising constant and $K_\nu(x)$ is a modified Bessel function of the second kind (defined in Appendix \ref{appa}).  Other parametrisations are given in \cite{gaunt vg}, \cite{mcc98}, and Section 4.1 of the book \cite{kkp01}, in which they use the name generalized Laplace distribution; these parametrisations are collected in a recent review \cite{vgsurvey} of the VG distribution. Alternative names include the McKay Type II distribution \cite{ha04} and Bessel function distribution \cite{m32}. The VG distribution is widely used in financial modelling \cite{mcc98,madan}, and has received recent interest as a natural limit distribution in probability theory (see, for example, \cite{aet21,gaunt vg}).
Further distributional properties and application areas are given in
\cite{vgsurvey,kkp01}. Henceforth, we set $\mu=0$.


The semi-heavy tails of the VG distribution make it suitable for financial modelling, and it is well-known that moments of all order exist. Lower order moments are readily obtained via the moment generating function; see \cite{mcc98} for the first four central moments. Moreover, a simple formula for the absolute moments of general order is available when $\beta=0$. For $Y\sim\mathrm{VG}(\nu,\alpha,0,0)$,
\begin{equation}\label{absmom}\mathbb{E}[|Y|^k]=\frac{2^{k}}{\sqrt{\pi}\alpha^k}\frac{\Gamma(\nu+(k+1)/2)\Gamma((k+1)/2)}{\Gamma(\nu+1/2)}, \quad k>k_*, 
\end{equation}
where $k_*=k_*(\nu):=\mathrm{max}\{-1,-2\nu-1\}$. This formula is readily obtained from the representation $Y=_d \alpha^{-1}Z\sqrt{G}$, where $Z\sim N(0,1)$ and $G\sim \Gamma(\nu+1/2,1/2)$ (with PDF $p_G(x)=x^{\nu-1/2}\mathrm{e}^{-x/2}/(2^{\nu+1/2}\Gamma(\nu+1/2))$, $x>0$) are independent (see \cite[Proposition 4.1.2]{kkp01}). Formula (\ref{absmom}) is a specialisation of our general formula (\ref{mom1}) for the case $\beta=0$; \cite[Proposition 2.2]{gaunt vg2} also gave a formula for $\mathbb{E}[|Y|^k]$, for $k>0$, with an unfortunate typo that our formula corrects.
When $k$ is an even integer, (\ref{absmom}) of course gives formulas for the moments of even order, whilst the moments of odd order are equal to zero on account of the symmetry of the distribution about the origin. However, in the general case $\beta\not=0$, the formulas available in the literature take a more complicated form. An exact formula for the moments of the VG distribution in terms of a finite sum involving gamma functions is given by \cite{scott}; this allows lower order moments to be calculated efficiently.  
A formula, involving a difference of two hypergeometric functions, for the moments of general order of the $\mathrm{VG}(\nu,\alpha,\beta,0)$ distribution was obtained by \cite[Theorem 2]{ha04} (they used the terminology McKay Type II distribution). As hypergeometric functions can be accurately and efficiently evaluated using modern computational algebra packages and mathematical software including the GNU Scientific Library, 
this is to date the simplest to implement and most practically 
useful formula
 in the literature for computing higher order moments.

In this paper, we obtain new formulas for the moments and absolute moments of the $\mathrm{VG}(\nu,\alpha,\beta,0)$ distribution in terms of a single hypergeometric function. We believe that these are the simplest possible formulas that can be given,
as the hypergeometric function does not reduce to a simpler form for our general parameter values and general values of the order of the raw moments $k>k_*$. On account of the wide use of the VG distribution in financial modelling and other application areas, it is useful for researchers to have these formulas at their disposal. To illustrate the convenience of our general formula, we efficiently derive the first few terms in the series expansions of the moments and absolute moments for the case that $|\beta/\alpha|\ll1$, a parameter regime often encountered when using the VG distribution to model log returns of financial assets (see, for instance, Example 3 of \cite{s04} concerning readings from the Dow-Jones Industrial Average for which $\beta/\alpha=0.0313$). We also deduce a formula for the moments and absolute moments of the product of two correlated zero mean normal random variables, and more generally the sum of independent copies of such random variables. These distributions themselves have numerous applications that date back to 1936 with the work of \cite{craig}; see \cite{gaunt22} for an overview of application areas and distributional properties.



\section{Results and proofs}

The following theorem is the main result of this paper. 
The formulas in this theorem involve the hypergeometric function, which is defined in Appendix \ref{appa}. 

\begin{theorem}\label{thm1} Let $X\sim \mathrm{VG}(\nu,\alpha,\beta,0)$, where $\nu > -1/2$, $0\leq|\beta|<\alpha$. Then, for $k>k_*$,
\begin{align}\label{mom1}\mathbb{E}[|X|^k]=\frac{2^k(1-\beta^2/\alpha^2)^{\nu+1/2}}{\sqrt{\pi}\alpha^{k}\Gamma(\nu+1/2)}\Gamma\Big(\nu+\frac{k+1}{2}\Big)\Gamma\Big(\frac{k+1}{2}\Big)\,{}_2F_1\bigg(\frac{k+1}{2},\nu+\frac{k+1}{2};\frac{1}{2};\frac{\beta^2}{\alpha^2}\bigg),
\end{align}
and, for odd $k\in\mathbb{Z}^+$,
\begin{align}\label{mom2}\mathbb{E}[X^k]=\frac{2^{k+1}\beta(1-\beta^2/\alpha^2)^{\nu+1/2}}{\sqrt{\pi}\alpha^{k+1}\Gamma(\nu+1/2)}\Gamma\Big(\nu+\frac{k}{2}+1\Big)\Gamma\Big(\frac{k}{2}+1\Big)\,{}_2F_1\bigg(\frac{k}{2}+1,\nu+\frac{k}{2}+1;\frac{3}{2};\frac{\beta^2}{\alpha^2}\bigg).
\end{align}
\end{theorem}

\begin{proof}
We first prove formula (\ref{mom1}). Let $k>k_*$. Using the formula (\ref{vgpdf}) for the VG PDF, the standard power series expansion of the exponential function, and interchanging the order of summation and integration, we obtain
\begin{align*}\mathbb{E}[|X|^k]&=M\int_{-\infty}^\infty \mathrm{e}^{\beta t}|t|^{\nu+k}K_{\nu}(\alpha|t|)\,\mathrm{d}t=M\sum_{i=0}^\infty\frac{\beta^i}{i!}\int_{-\infty}^\infty t^i|t|^{\nu+k}K_\nu(\alpha|t|)\,\mathrm{d}t\\
&=M\sum_{j=0}^\infty\frac{\beta^{2j}}{(2j)!}\int_{-\infty}^\infty |t|^{\nu+k+2j}K_\nu(\alpha|t|)\,\mathrm{d}t=2M\sum_{j=0}^\infty\frac{\beta^{2j}}{(2j)!}\int_0^\infty t^{\nu+k+2j}K_\nu(\alpha t)\,\mathrm{d}t.
\end{align*}
Note that the expectation $\mathbb{E}[|X|^k]$ exists for $k>k_*$ due to the limiting forms (\ref{Ktend0}) and (\ref{Ktendinfinity}).
Evaluating the final integral using (\ref{integral}) now yields the formula
\begin{align}\mathbb{E}[|X|^k]&=2M\sum_{j=0}^\infty\frac{(\beta^2)^j}{(2j)!}\cdot \frac{2^{\nu+k+2j-1}}{\alpha^{\nu+k+2j+1}}\Gamma\Big(\nu+j+\frac{k+1}{2}\Big)\Gamma\Big(j+\frac{k+1}{2}\Big)\nonumber\\
&=\frac{2^k(1-\beta^2/\alpha^2)^{\nu+1/2}}{\sqrt{\pi}\alpha^{k}\Gamma(\nu+1/2)}\Gamma\Big(\nu+\frac{k+1}{2}\Big)\Gamma\Big(\frac{k+1}{2}\Big)\times\nonumber\\
\label{momint}&\quad\times\sum_{j=0}^\infty \frac{((k+1)/2)_j(\nu+(k+1)/2)_j}{(1/2)_j}\frac{(\beta^2/\alpha^2)^j}{j!},
\end{align}
where in obtaining the second equality we used the standard formula $(u)_k=\Gamma(u+k)/\Gamma(u)$ to put the infinite series into the hypergeometric form (\ref{gauss}). Formula (\ref{mom1}) now follows from (\ref{momint}) and the definition of the hypergeometric function (\ref{gauss}).

Now suppose $k$ is an odd integer. The proof of formula (\ref{mom2}) is similar to that of (\ref{mom1}), so we provide a more efficient exposition. By a similar reasoning to before we obtain
\begin{align*}\mathbb{E}[X^k]&=M\sum_{i=0}^\infty\frac{\beta^i}{i!}\int_{-\infty}^\infty t^{k+i}|t|^{\nu}K_\nu(\alpha|t|)\,\mathrm{d}t=2M\sum_{j=0}^\infty\frac{\beta^{2j+1}}{(2j+1)!}\int_0^\infty t^{\nu+k+2j+1}K_\nu(\alpha t)\,\mathrm{d}t.
\end{align*}
Formula (\ref{mom2}) is now readily obtained by evaluating the integral using (\ref{Ktendinfinity}) and putting the resulting infinite series into the hypergeometric form. We omit the details. \end{proof}

\begin{remark}Let $k\in\mathbb{Z}^+$, and define $\ell:=\lceil k/2\rceil+1/2$ and $m:=k\,\mathrm{mod}\,2$. Then, on combining (\ref{mom1}) and (\ref{mom2}) we obtain the following compact formula for the moments of the $ \mathrm{VG}(\nu,\alpha,\beta,0)$ distribution:
\begin{equation}\label{mom3}\mathbb{E}[X^k]=\frac{2^{k}(2\beta/\alpha)^m(1-\beta^2/\alpha^2)^{\nu+1/2}}{\sqrt{\pi}\alpha^{k}\Gamma(\nu+1/2)}\Gamma(\nu+\ell)\Gamma(\ell)\,{}_2F_1\bigg(\ell,\nu+\ell;\frac{1}{2}+m;\frac{\beta^2}{\alpha^2}\bigg).
\end{equation}
\end{remark}

\begin{remark}As mentioned in the Introduction, the parameter regime $|\beta/\alpha|\ll1$ arises in financial modelling. Applying the series representation (\ref{gauss}) of the hypergeometric function and the generalised binomial series of $(1-\beta^2/\alpha^2)^{\nu+1/2}$ to (\ref{mom1}) and (\ref{mom3}), we can obtain series expansions for the moments and absolute moments of the $ \mathrm{VG}(\nu,\alpha,\beta,0)$ distribution in this parameter regime. We report the first few terms. For $k>k_*$,
\begin{align*}\mathbb{E}[|X|^k]&=\frac{(2/\alpha)^k}{\sqrt{\pi}\Gamma(\nu+1/2)}\Gamma\Big(\nu+\frac{k+1}{2}\Big)\Gamma\Big(\frac{k+1}{2}\Big)\bigg[1+\frac{k(k+2\nu+2)}{2}\frac{\beta^2}{\alpha^2}\\
&\quad+\frac{1}{24}k(k^3+4k^2(\nu+2)+4k(\nu^2+3\nu+4)-8\nu^2+8\nu+12)\frac{\beta^4}{\alpha^4}+O\bigg(\frac{\beta^6}{\alpha^6}\bigg)\bigg],
\end{align*}
as $\beta/\alpha\rightarrow0$, whilst, for $k\in\mathbb{Z}^+$,
\begin{align*}\mathbb{E}[X^k]&=\frac{(2/\alpha)^{k}(2\beta/\alpha)^m}{\sqrt{\pi}\Gamma(\nu+1/2)}\Gamma(\nu+\ell)\Gamma(\ell)\bigg[1+\bigg(\frac{2\ell(\ell+\nu)}{2m+1}-\nu-\frac{1}{2}\bigg)\frac{\beta^2}{\alpha^2}\\
&\quad+\bigg(\frac{2\ell(\ell+1)(\ell+\nu)(\ell+\nu+1)}{4m^2+8m+3}-\frac{\ell(2\nu+1)(\ell+\nu)}{2m+1}+\frac{4\nu^2-1}{8}\bigg)\frac{\beta^4}{\alpha^4}+O\bigg(\frac{\beta^6}{\alpha^6}\bigg)\bigg],
\end{align*}
as $\beta/\alpha\rightarrow0$.
\end{remark}

Now, let $(U,V)$ be a bivariate normal random vector with  zero mean vector, variances $(\sigma_U^2,\sigma_V^2)$ and correlation coefficient $\rho$. Set $s=\sigma_U\sigma_V$, and let $Z=UV$ be the product of these correlated normal random variables.  Consider also the mean $\overline{Z}_n=n^{-1}(Z_1+Z_2+\cdots+Z_n)$, where $Z_1,Z_2,\ldots,Z_n$ are independent copies of $Z$. It was noted by \cite{gaunt thesis} that $Z$ has a VG distribution and later by \cite{gaunt prod} that more generally
\begin{equation}\label{vgrep}\overline{Z}_n\sim\mathrm{VG}\bigg(\frac{n-1}{2},\frac{n}{s(1-\rho^2)},\frac{n\rho}{s(1-\rho^2)},0\bigg).
\end{equation}
 Combining (\ref{vgrep}) with (\ref{mom1}) and (\ref{mom3}), the following formulas for the moments and absolute moments of $\overline{Z}_n$ are immediate; formulas for the moments of $Z$ follow on setting $n=1$. These formulas are simpler than those recently given by \cite{gaunt22}.

\begin{corollary}Let the previous notations prevail. Then, for $k>-1$,
\begin{equation*}\mathbb{E}[|\overline{Z}_n|^k]=\frac{(2s/n)^k(1-\rho^2)^{n/2+k}}{\sqrt{\pi}\Gamma(n/2)}\Gamma\Big(\frac{n+k}{2}\Big)\Gamma\Big(\frac{k+1}{2}\Big)\,{}_2F_1\bigg(\frac{k+1}{2},\frac{n+k}{2};\frac{1}{2};\rho^2\bigg),
\end{equation*}
whilst, for $k\in\mathbb{Z}^+$,
\begin{equation*}\mathbb{E}[\overline{Z}_n^k]=\frac{(2s/n)^k(2\rho)^m(1-\rho^2)^{n/2+k}}{\sqrt{\pi}\Gamma(n/2)}\Gamma\Big(\frac{n-1}{2}+\ell\Big)\Gamma(\ell)\,{}_2F_1\bigg(\ell,\frac{n-1}{2}+\ell;\frac{1}{2}+m;\rho^2\bigg).
\end{equation*}
\end{corollary}



\appendix

\section{Special functions}\label{appa}
In this appendix, we define the modified Bessel function of the second kind and the hypergeometric function, and present some basic properties that are used in this paper. All properties can be found in \cite{olver}. The modified Bessel function of the second kind $K_\nu(x)$ is defined, for $\nu\in\mathbb{R}$ and $x>0$, by
\[K_\nu(x)=\int_0^\infty \mathrm{e}^{-x\cosh(t)}\cosh(\nu t)\,\mathrm{d}t.
\]
It has the following asymptotic behaviour:
\begin{eqnarray}\label{Ktend0}K_{\nu} (x) &\sim& \begin{cases} 2^{|\nu| -1} \Gamma (|\nu|) x^{-|\nu|}, & \quad x \downarrow 0, \: \nu \not= 0, \\
-\log x, & \quad x \downarrow 0, \: \nu = 0, \end{cases} \\
\label{Ktendinfinity} K_{\nu} (x) &\sim& \sqrt{\frac{\pi}{2x}} \mathrm{e}^{-x}, \quad x \rightarrow \infty,\: \nu\in\mathbb{R}.
\end{eqnarray}
For $r>|\nu|$, we have the following definite integral formula:
\begin{equation}\label{integral}\int_0^\infty t^{r-1}K_\nu(t)\,\mathrm{d}t=2^{r-2}\Gamma\Big(\frac{r-\nu}{2}\Big)\Gamma\Big(\frac{r+\nu}{2}\Big).
\end{equation}
The (Gaussian) hypergeometric function is defined, for $|x|<1$, by the power series
\begin{equation}\label{gauss}{}_2F_1(a,b;c;x)=\sum_{j=0}^\infty\frac{(a)_j(b)_j}{(c)_j}\frac{x^j}{j!},
\end{equation}
where $(u)_j=u(u+1)\cdots(u+k-1)$ is the ascending factorial. 

\section*{Acknowledgements}
I would like to thank the reviewers for carefully reading the manuscript and for their helpful comments and suggestions.

\end{document}